\documentclass[10pt]{amsart}
\usepackage[all]{xy}
\usepackage{amsmath,amsthm,amsfonts,amssymb}
\newtheorem{theorem}{Theorem}[section]
\newtheorem{proposition}[theorem]{Proposition}
\newtheorem{corollary}[theorem]{Corollary}
\newtheorem{lemma}[theorem]{Lemma}
\newtheorem{remark}[theorem]{Remark}
\newtheorem{example}[theorem]{Example}
\newcommand{\A}{{\mathbb A}}
\newcommand{\C}{{\mathcal C}}
\newcommand{\Q}{{\mathbb Q}}
\newcommand{\ST}{{\mathcal S}}
\newcommand{\T}{{\mathcal T}}
\newcommand{\Z}{{\mathbb Z}}
\newcommand{\LL}{{\mathbb L}}
\newcommand{\LLdn}{{\mathbb L}_{\bullet}}
\newcommand{\LLup}{{\mathbb L}^{\bullet}}
\newcommand{\di}{\displaystyle}

\title{A composition formula for manifold structures}
\author{Andrew Ranicki}
\address{School of Mathematics\newline
\indent University of Edinburgh\newline \indent Edinburgh EH9
3JZ\newline \indent Scotland, UK} \email{\noindent
A.Ranicki@ed.ac.uk}



\date{\today}

\begin{document}
\begin{abstract}
The structure set $\ST^{TOP}(M)$ of an $n$-dimensional topological
manifold $M$ for $n \geqslant 5$ has a homotopy invariant functorial
abelian group structure, by the algebraic version of the
Browder-Novikov-Sullivan-Wall surgery theory.
An element $(N,f) \in \ST^{TOP}(M)$ is an equivalence class of
$n$-dimensional manifolds $N$ with a homotopy equivalence $f:N \to M$.
The composition formula is that $(P,fg)=(N,f)+f_*(P,g) \in \ST^{TOP}(M)$
for homotopy equivalences $g:P \to N$, $f:N \to M$. The formula is required
for a paper of Kreck and L\"uck.
\end{abstract}
\maketitle \hfill {\it For Prof. F. Hirzebruch}

\section*{Introduction}

The structure set $\ST^{TOP}(M)$ of an $n$-dimensional topological
manifold $M$ is the pointed set of equivalence classes of
pairs $(N,f)$ with $N$ an $n$-dimensional manifold and
$f:N \to M$ a homotopy equivalence, with
$(N_1,f_1)=(N_2,f_2) \in
\ST^{TOP}(M)$ if and only if $(f_1)^{-1}f_2:N_2 \to N_1$ is homotopic
to a homeomorphism, and $(M,1) \in \ST^{TOP}(M)$ the base point.
The equivalence class $(N,f) \in \ST^{TOP}(M)$
is called the {\it structure invariant} of $(N,f)$.
For $n \geqslant 5$ $\ST^{TOP}(M)$ is determined by the fundamental
exact sequence of pointed sets of the Browder-Novikov-Sullivan-Wall surgery
theory
$$\xymatrix@C-5pt{
\dots \ar[r] & L_{n+1}(\Z[\pi_1(M)]) \ar[r] & \ST^{TOP}(M) \ar[r]^-{\di{\eta}} &
\T^{TOP}(M) \ar[r]^-{\di{\theta}} & L_n(\Z[\pi_1(M)])}$$
with $\T^{TOP}(M)=[M,G/TOP]$ the normal bordism set,
$G/TOP$ the classifying space for fibre homotopy trivialized
topological bundles, $\eta$ sending a structure invariant to its
bordism class when regarded as a normal map, and $\theta$ the
surgery obstruction function.

This paper will only consider homotopy equivalences of manifolds which are
simple, i.e. have zero Whitehead torsion, in order to take advantage of the
$s$-cobordism theorem. However, all the results have analogues for
homotopy equivalences which are not necessarily simple.

The algebraic surgery sequence of abelian groups
$$\xymatrix@C-5pt{
\dots \ar[r] & L_{n+1}(\Z[\pi_1(M)]) \ar[r] & \ST_{n+1}(M) \ar[r] &
H_n(M;\LLdn) \ar[r]^-{\di{A}} & L_n(\Z[\pi_1(M)])}$$
was constructed in Ranicki \cite[\S16]{ranicki4} for any space $M$,
with $\LLdn$ the 1-connective quadratic $L$-theory spectrum of $\Z$
such that $\LL_0 \simeq G/TOP$, and with $A$ the assembly map.
The algebraic and topological surgery sequences are isomorphic,
with bijections
$$\begin{array}{l}
\ST^{TOP}(M) \to \ST_{n+1}(M)~;~(N,f) \mapsto s(f)~,\\[1ex]
\T^{TOP}(M)=[M,G/TOP] \to H_n(M;\LLdn)~.
\end{array}$$

\S1 recalls the algebraic and topological surgery exact sequences.

\S2 uses the algebraic surgery exact sequence to prove the composition formula
$$(P,fg)~=~(N,f)+f_*(P,g) \in \ST^{TOP}(M)\eqno{(*)}$$
for homotopy equivalences $f:N \to M$, $g:P \to N$
of $n$-dimensional topological manifolds (Theorem \ref{composite}).
The addition in the structure set is thus given by
$$(N,f)+(N',f')~=~(P,fg) \in \ST^{TOP}(M)$$
with $(P,g)=(f_*)^{-1}(N',f') \in \ST^{TOP}(N)$.

A homotopy equivalence of $n$-dimensional manifolds $f:N \to M$
is automatically a normal map $(f,b)$, with
$$b~:~\nu_N \to (f^{-1})^*\nu_N~=~\nu_M\oplus \nu_f~.$$
Use the normal invariant $(f,b)\in [M,G/TOP]$ to define a function
$$\eta~:~\ST^{TOP}(M)~=~\ST_{n+1}(M) \to \T^{TOP}(M)~;~
(N,f)=s(f) \mapsto \eta(N,f)=(f,b)~.$$
The composite
of homotopy equivalences $f:N \to M$, $g:P \to N$ is a homotopy
equivalence $fg:P \to M$ with normal invariant
$$\eta(P,fg)~=~\eta(N,f) \oplus (f^{-1})^*\eta(P,g) \in \T^{TOP}(M)\eqno{(**)}$$
by the sum formula of Brumfiel \cite{brumfiel} and
Madsen-Taylor-Williams \cite{mtw}.

At first sight, there is a disparity between the sum formulae $(*)$ and
$(**)$, since it is well-known that the disjoint union
addition + and the Whitney sum addition $\oplus$ on the surgery
classifying space $G/TOP \simeq \LL_0$ are different, and the normal
bordism set $\T^{TOP}(M)=[M,G/TOP]$ has two abelian group structures.
A map $f:N \to M$ of $n$-dimensional manifolds induces a function
$f_*:\T^{TOP}(N) \to \T^{TOP}(M)$ which is a morphism of abelian groups
with respect to $+$, and a function $f^*:\T^{TOP}(M) \to \T^{TOP}(N)$ which is
a morphism of abelian groups with respect to $\oplus$.  If $f$ is a
homotopy equivalence then both $f_*$ and $f^*$ are isomorphisms, but in
general they are not inverse bijections.  The diagram
$$\xymatrix@C+10pt{\ST^{TOP}(N)\ar[d]_-{\di{\eta}}
\ar[r]^-{\di{f_*}}_-{\di{\cong}} & \ST^{TOP}(M) \ar[d]^-{\di{\eta}} \\
\T^{TOP}(N) \ar[r]^-{\di{(f^{-1})^*}}_-{\di{\cong}} & \T^{TOP}(M)}$$
is not commutative, essentially since $f$ does not preserve the
Hirzebruch ${\mathcal L}$-genus ${\mathcal L}(M) \in H_{n-4*}(M;\Q)$
$$f_*{\mathcal L}(N)\neq {\mathcal L}(M) \in H_{n-4*}(M;\Q)~.$$
An $n$-dimensional manifold $M$ has a canonical
$\LLup$-theory fundamental class $[M]_{\LL} \in H_n(M;\LLup)$ with $\LLup$ the
symmetric $L$-theory spectrum of $\Z$ (Ranicki \cite[\S16]{ranicki4}), such that
$$[M]_{\LL}\otimes \Q~=~{\mathcal L}(M) \in
H_n(M;\LLup)\otimes \Q~=~H_{n-4*}(M;\Q)~.$$
The failure of commutativity in the diagram
is traced in \S3 to the failure (in general) of $f$ to preserve
the symmetric $L$-theory fundamental classes, with
$$f_*[N]_{\LL}-[M]_{\LL}~=~(1+T)\eta(N,f) \in
{\rm im}(1+T:H_n(M;\LLdn) \to H_n(M;\LLup))~.$$
The normal invariant function $\eta:\ST^{TOP}(M) \to \T^{TOP}(M)$
is a morphism of abelian groups with respect to +, so $(*)$ gives
$$\eta(P,fg)~=~\eta(N,f)+\eta(f_*(P,g)) \in \T^{TOP}(M)~.$$
The symmetric $L$-theory orientation of manifolds is used in
Corollary \ref{reconcile} to prove that
$$\eta(N,f)+\eta (f_*(P,g))~=~\eta(N,f)\oplus (f^{-1})^*\eta(P,g)
\in \T^{TOP}(M)$$
reconciling $(*)$ and $(**)$.

This paper is a substantial expansion of the paper posted on the arXiv
in November 2006 (math.AT/0608705). I am grateful to Larry Taylor for
asking me about the relationship between the manifold structure
composition formula (*) and the normal invariant composition formula (**)
which stimulated me to expand the original version of the paper.

\section{The algebraic and topological surgery exact sequences}

We refer to Wall \cite{wall1} for the basic properties of geometric
Poincar\'e complexes and pairs.  A {\it Spivak normal structure}
$(\nu_M,\rho_M)$ on an $n$-dimensional geometric Poincar\'e complex $M$
is a $(k-1)$-spherical fibration $\nu_M:M \to BG(k)$ together with a
map $\rho_M:S^{n+k} \to T(\nu_M)$ to the Thom space with
Hurewicz-Thom image the fundamental class
$$h(\rho_M)~=~[M] \in \dot H_{n+k}(T(\nu_M))~=~H_n(M)~.$$
Again, we refer to \cite{wall1} for the existence and uniqueness
properties of Spivak normal structures.  An $n$-dimensional manifold
$M$ is an $n$-dimensional geometric Poincar\'e complex.  An embedding
$M \subset S^{n+k}$ ($k$ large) has a topological normal bundle
$\nu_M:M \to BTOP(k)$ and the Pontrjagin-Thom construction gives a map
$\rho_M:S^{n+k}\to T(\nu_M)$ such that $(J\nu_M,\rho_M)$ is a Spivak
normal structure, with $J:BTOP(k) \to BG(k)$ the forgetful map.
Similarly for geometric Poincar\'e pairs and manifolds with boundary.

A {\it normal map} $(f,b):N \to M$ of $n$-dimensional geometric
Poincar\'e complexes with Spivak normal structures $(\nu_N,\rho_N)$,
$(\nu_M,\rho_M)$ is a degree 1 map $f:N \to M$ together with a map
$b:\nu_N \to \nu_M$ of Spivak normal fibrations such that
$$T(b)\rho_N~\simeq~\rho_M~:~S^{n+k} \to T(\nu_M)~.$$
\indent
A {\it topological normal map} $(f,b):N \to M$ from an $n$-dimensional
manifold $N$ to an $n$-dimensional geometric Poincar\'e complex $M$ is
a degree 1 map $f:N \to M$ together with a bundle map $b:\nu_N \to
\nu_M\oplus \nu_b$ for some topological bundle $\nu_b:M \to BTOP(k)$ with
a fibre homotopy trivialization $J\nu_b \simeq *:M \to BG(k)$, in which
case $(f,Jb):N \to M$ is a normal map of the underlying geometric
Poincar\'e complexes with $Jb:J\nu_N \to J\nu_M$.
In particular, a homotopy equivalence $f:N \to M$ of
manifolds determines a topological normal map $(f,b)$ with
$b:\nu_N \to (f^{-1})^*\nu_N$.

As usual, let
$$G/TOP~=~\text{homotopy fibre of}~J:BTOP \to BG~,$$
the classifying space for fibre homotopy trivialized stable topological
bundles.  A topological normal map $(f,b):N \to M$ has a
{\it topological normal invariant} $\eta(f,b) \in [M,G/TOP]$ such that
$b:\nu_M \to \nu_N \oplus \nu_b$ with
$$\nu_b~:~M \xymatrix{\ar[r]^-{\di{\eta(f,b)}}&} G/TOP \to BTOP~.$$
The {\it normal map set} $\T^{TOP}(M)$ of a manifold $M$
is the pointed set of bordism classes of topological
normal maps $(f,b):N \to M$. The function
$$\T^{TOP}(M) \to [M,G/TOP]~;~(f,b) \mapsto \eta(f,b)$$
is a bijection. The inverse is given by the Browder-Novikov transversality
construction: given a map $\eta:M \to G/TOP$ let
$\nu_b:M \to BTOP$ be the corresponding fibre homotopy trivial topological
bundle, let  $h_\eta:T(\nu_M) \simeq
T(\nu_M\oplus \nu_b)$ be the corresponding homotopy equivalence of
Thom spaces, and make the map
$$h_{\eta}\rho_M~:~S^{n+k} \to T(\nu_M)~\simeq~T(\nu_M\oplus \nu_b)$$
topologically transverse regular at the zero section $M \subset T(\nu_M\oplus \nu_b)$,
to obtain a topological normal map
$$(f,b)~=~h_{\eta}\rho_M\vert~:~(N,\nu_N)~=~
((h_{\eta}\rho_M)^{-1}(M),\nu_{N \subset S^{n+k}}) \to (M,\nu_M\oplus \nu_b)$$
with topological normal invariant $\eta$.
For a manifold with boundary $(M,\partial M)$ there is a relative normal map
set $\T(M,\partial M)$ of topological normal maps
$(f,b):(N,\partial N) \to (M,\partial M)$  with a bijection
$$\T^{TOP}(M,\partial M) ~\cong~ [M,G/TOP]~.$$
There is also a rel $\partial$ normal map
set $\T^{TOP}_{\partial}(M,\partial M)$ of topological normal maps
$(f,b):N \to M$ such that $\partial f:\partial N \to \partial M$ is a
homeomorphism, with a bijection
$$\T^{TOP}_{\partial}(M,\partial M) ~\cong~ [M,\partial M;G/TOP,*]~.$$
\indent
The Wall \cite{wall2} surgery obstruction groups $L_n(\Lambda)$ of a ring
$\Lambda$ with involution $\lambda \mapsto \overline{\lambda}$
were expressed in
Ranicki \cite{ranicki2} as the cobordism groups of $n$-dimensional
quadratic Poincar\'e complexes $(C,\psi)$ over $\Lambda$, with $C$
an $n$-dimensional based f.g. free $\Lambda$-module chain complex.
The quadratic structure is an element
$$\psi \in Q_n(C)~=~H_n(W\otimes_{\Z[\Z_2]}(C\otimes_{\Lambda}C))$$
where
$$\xymatrix{W~:~\dots \ar[r] & \Z[\Z_2] \ar[r]^-{1-T} &
\Z[\Z_2] \ar[r]^-{1+T} &\Z[\Z_2] \ar[r]^-{1-T} &\Z[\Z_2]}$$ is the
standard free $\Z[\Z_2]$-module resolution of $\Z$, and $T \in \Z_2$
acts on
$$C\otimes_{\Lambda}C~=~C\otimes_{\Z}C/
\{x\otimes \lambda y - \overline{\lambda}x \otimes y\,\vert\,
x,y \in C,\,\lambda \in \Lambda\}$$
by
$$T~:~C_p\otimes_{\Lambda}C_q \to C_q\otimes_{\Lambda}C_p~;~
x \otimes y \mapsto (-)^{pq} y \otimes x~.$$
The $\Lambda$-module chain map
$$(1+T)\psi_0~:~C^{n-*}~=~{\rm Hom}_{\Lambda}(C,\Lambda)_{*-n} \to C$$
is required to be a simple chain equivalence. The quadratic
$L$-groups are the homotopy groups of an $\Omega$-spectrum
$\LLdn(\Lambda)$ of quadratic Poincar\'e complexes over $\Lambda$, with
$$\pi_n(\LLdn(\Lambda))~=~L_n(\Lambda)~.$$
For the surgery exact sequence it is necessary to consider the
1-connective cover of $\LLdn(\Z)$, which is denoted by $\LLdn$,
and is such that
$${\mathbb L}_0~\simeq~G/TOP~,~
\pi_n(\LLdn)~=~\pi_n(G/TOP)~=~
L_n(\Z)~=~\begin{cases} \Z &{\rm if}~n \equiv 0(\bmod\,4)\\
0 &{\rm if}~n \equiv 1(\bmod\,4)\\
\Z_2 &{\rm if}~n \equiv 2(\bmod\,4)\\
0 &{\rm if}~n \equiv 3(\bmod\,4)~.
\end{cases}$$
The $H$-space structure on ${\mathbb L}_0$ given by the direct sum of quadratic
Poincar\'e complexes corresponds to the disjoint union $H$-space
structure on $G/TOP$.

The {\it symmetric $L$-groups} $L^n(\Lambda)$ of a ring $\Lambda$ with
involution were defined in \cite{ranicki2} (following Mishchenko) as
the cobordism groups of $n$-dimensional symmetric Poincar\'e complexes
$(C,\phi)$ over $\Lambda$, with $C$ an $n$-dimensional based f.g.  free
$\Lambda$-module chain complex. The symmetric structure is an element
$$\phi \in Q^n(C)~=~H_n({\rm Hom}_{\Z[\Z_2]}(W,C\otimes_{\Lambda}C))$$
and the $\Lambda$-module chain map $\phi_0:~C^{n-*} \to C$
is required to be a simple chain equivalence. The symmetrization maps
$$1+T~:~L_n(\Lambda) \to L^n(\Lambda)~;~(C,\psi) \mapsto (C,(1+T)\psi)$$
are isomorphisms modulo 8-torsion (\cite[Proposition 8.2]{ranicki1}).
For any rings with involution $\Lambda,\Lambda'$
the tensor product over $\Z$ defines $L$-theory products
$$\begin{array}{l}
L^n(\Lambda) \otimes_{\Z} L^{n'}(\Lambda') \to L^{n+n'}(\Lambda\otimes_{\Z}\Lambda')~,\\[1ex]
L^n(\Lambda) \otimes_{\Z} L_{n'}(\Lambda') \to L_{n+n'}(\Lambda\otimes_{\Z}\Lambda')
\end{array}$$
and also
$$L_n(\Lambda) \otimes_{\Z} L_{n'}(\Lambda') \to L_{n+n'}(\Lambda\otimes_{\Z}\Lambda')~;~
x \otimes x' \mapsto (1+T)x\otimes x' = x \otimes (1+T)x'~.$$
\indent The {\it symmetric signature} of an $n$-dimensional geometric
Poincar\'e complex $M$ was defined in Ranicki \cite{ranicki3} (again,
following Mishchenko)  to be the cobordism class
$$\sigma^*(M)~=~(C(\widetilde{M}),\phi_M) \in L^n(\Z[\pi_1(M)])$$
of the $n$-dimensional symmetric Poincar\'e complex
$(C(\widetilde{M}),\phi_M)$ over $\Z[\pi_1(M)]$ with $\widetilde{M}$ is
the universal cover of $M$.
The symmetric $L$-groups are the homotopy groups of an $\Omega$-spectrum
$\LLup(\Lambda)$ of symmetric Poincar\'e complexes over $\Lambda$, with
$$\pi_n(\LLup(\Lambda))~=~L^n(\Lambda)~.$$
In particular, for $\Lambda=\Z$ we write $\LLup(\Z)=\LLup$, with
$$\pi_n(\LLup)~=~L^n(\Z)~=~
\begin{cases} \Z &{\rm if}~n \equiv 0(\bmod\,4)\\
\Z_2 &{\rm if}~n \equiv 1(\bmod\,4)\\
0 &{\rm if}~n \equiv 2(\bmod\,4)\\
0 &{\rm if}~n \equiv 3(\bmod\,4)~.
\end{cases}$$
The tensor product of symmetric Poincar\'e complexes over $\Z$
gives $\LLup$ the structure of a ring spectrum (with a 1), $\LLdn(\Lambda)$ is an
$\LLup$-module spectrum, and $\LLdn$ is also a ring spectrum (but only with an 8, not a 1).
A topological bundle has an
$\LLup$-orientation, and an $n$-dimensional manifold $M$ has a
symmetric $L$-theory fundamental class $[M]_{\LL} \in
H_n(M;\LLup)$ (\cite[Proposition 16.16]{ranicki4}) with assembly
$$A([M]_{\LL})~=~\sigma^*(M) \in L^n(\Z[\pi_1(M)])~.$$
\indent The {\it quadratic signature} of a
normal map $(f,b):N \to M$ of $n$-dimensional geometric Poincar\'e complexes
was defined in Ranicki \cite{ranicki3} to be the cobordism class
$$\sigma_*(f,b)~=~(\C(f^!),\psi_b) \in L_n(\Z[\pi_1(M)])$$
of the {\it quadratic kernel} $n$-dimensional quadratic Poincar\'e complex
$(\C(f^!),\psi_b)$ over $\Z[\pi_1(M)]$, constructed as follows.
The Umkehr $\Z[\pi_1(M)]$-module chain map is the composite
$$f^!~:~C(\widetilde{M}) \simeq C(\widetilde{M})^{n-*}
\xymatrix{\ar[r]^-{\di{\widetilde{f}^*}} & }
C(\widetilde{N})^{n-*} \simeq C(\widetilde{N})~,$$
where $\widetilde{M}$ is the universal cover of $M$,
$\widetilde{N}=f^*\widetilde{M}$ is the pullback cover of $N$, and
$\widetilde{f}:\widetilde{N} \to \widetilde{M}$ is a
$\pi_1(M)$-equivariant lift of $f$. The algebraic mapping cone
$\C(f^!)$ fits into the chain homotopy direct sum system
$$\xymatrix{C(\widetilde{M}) \ar@<1ex>[r]^-{\di{f^!}}
\ar@{<-}@<-1ex>[r]_{\di{\widetilde{f}}} & C(\widetilde{N})
 \ar@<1ex>[r]^-{\di{e}}\ar@{<-}@<-1ex>[r] & \C(f^!)}$$
with $e:C(\widetilde{N}) \to \C(f^!)$ the inclusion and
$\widetilde{f} f^!~\simeq~1~:~ C(\widetilde{M}) \to C(\widetilde{M})$.
Thus there is defined a chain equivalence of symmetric Poincar\'e
complexes
$$(C(\widetilde{N}),\phi_N)~\simeq~(\C(f^!),(1+T)\psi_b)
\oplus (C(\widetilde{M}),\phi_M)~,$$
and
$$\begin{array}{l}
H_*(\C(f^!))~=~K_*(N)~=~{\rm ker}(\widetilde{f}_*:H_*(\widetilde{N})
\to H_*(\widetilde{M}))~,\\[1ex]
H_*(\widetilde{N})~=~ H_*(\C(f^!)) \oplus H_*(\widetilde{M})~.
\end{array}$$
Let $\nu_{\widetilde{M}}$ be the pullback of $\nu_M$ to the universal
cover $\widetilde{M}$ of $M$, and let $\nu_{\widetilde{N}}$ be the
pullback of $\nu_N$ to the pullback cover $\widetilde{N}=f^*\widetilde{M}$
of $N$, so that $b$ lifts to a $\pi_1(M)$-equivariant map $\widetilde{b}:
\nu_{\widetilde{N}}\to\nu_{\widetilde{M}}$.
The $\pi_1(M)$-equivariant $S$-dual of
$T(\widetilde{b}):T(\nu_{\widetilde{N}})\to T(\nu_{\widetilde{M}})$
is a $\pi_1(M)$-equivariant stable map
$F:\Sigma^{\infty}\widetilde{M}^+\to\Sigma^{\infty}\widetilde{N}^+$
inducing $f^!$ on the chain level. The {\it quadratic construction} of
\cite[\S1]{ranicki3} is a natural transformation
$$\psi_F~:~C(M) \to W\otimes_{\Z[\Z_2]}(C(\widetilde{N})\otimes_{\Z[\pi_1(M)]}
C(\widetilde{N}))$$
and the quadratic structure in $\sigma_*(f,b)$ is defined by
$$\psi_b~=~(e\otimes e)\psi_F[M] \in Q_n(\C(f^!))~.$$
There is also a version for a normal map of geometric Poincar\'e pairs.

The surgery obstruction $\sigma_*(f,b)\in L_n(\Z[\pi_1(M)])$ was
defined by Wall \cite{wall2} for a topological  $n$-dimensional normal map
$(f,b):N \to M$ of manifolds with boundary such that
$\partial f:\partial N \to \partial M$ is a homotopy equivalence.
The surgery obstruction was identified in \cite{ranicki3} with
the quadratic signature of the underlying normal map of
geometric Poincar\'e pairs.

The {\it algebraic normal invariant} of a topological normal map
$(f,b):N \to M$ of $n$-dimensional manifolds is the cobordism class
$$t(f,b) ~=~(C,\psi) \in H_n(M;\LLdn)$$
defined in \cite[Proposition 18.3]{ranicki4} (as recalled in \S2 below),
with symmetrization
$$(1+T)t(f,b)~=~f_*[N]_{\LL}-[M]_{\LL} \in H_n(M;\LLdn)$$
and assembly the surgery obstruction
$$A(t(f,b))~=~\sigma_*(f,b) \in L_n(\Z[\pi_1(M)])~.$$
\indent
The {\it structure set} $\ST^{TOP}(M)$ of an
$n$-dimensional topological manifold $M$ was already defined
in the Introduction. For a manifold with boundary $(M,\partial M)$
there is a relative structure set $\ST^{TOP}(M,\partial M)$
of equivalence classes of
pairs $(N,f)$ with $N$ an $n$-dimensional manifold with boundary and
$f:N \to M$ a homotopy equivalence.
There is also a rel $\partial$ structure set $\ST^{TOP}_{\partial}(M,\partial M)$
of equivalence classes of
pairs $(N,f)$ with $N$ an $n$-dimensional manifold with boundary and
$f:N \to M$ a homotopy equivalence such that
$\partial f:\partial N \to \partial M$ is a homeomorphism.

For $n \geqslant 5$ and an $n$-dimensional manifold $M$ the 
Browder-Novikov-Sullivan-Wall surgery theory (extended to the
topological category by Kirby and Siebenmann \cite{ks}) fits
$\ST^{TOP}_{\partial}(M\times D^k)$ ($k \geqslant 0$) into
an exact sequence of pointed sets
$$\begin{array}{l}
\dots \xymatrix@C-2pt{\ar[r]&} \ST^{TOP}_{\partial}(M \times D^k)
\xymatrix@C-2pt{\ar[r]^-{\di{\eta}}&}
\T^{TOP}_{\partial}(M\times D^k) \xymatrix@C-2pt{\ar[r]^-{\di{\theta}}&}
L_{n+k}(\Z[\pi_1(M)]) \xymatrix@C-2pt{\ar[r]&}\\[1ex]
\hphantom{\dots} \dots \xymatrix@C-2pt{\ar[r]&}
\ST^{TOP}_{\partial}(M \times D^1) \xymatrix@C-2pt{\ar[r]^-{\di{\eta}}&}
\T^{TOP}_{\partial}(M\times D^1) \xymatrix@C-2pt{\ar[r]^-{\di{\theta}}&}
L_{n+1}(\Z[\pi_1(M)])\\[1ex]
\hphantom{\dots\dots} \xymatrix@C-2pt{\ar[r]&} \ST^{TOP}(M) 
\xymatrix@C-2pt{\ar[r]^-{\di{\eta}}&}
\T^{TOP}(M) \xymatrix@C-2pt{\ar[r]^-{\di{\theta}}&} L_n(\Z[\pi_1(M)])
\end{array}$$
which is in fact an exact sequence of abelian groups until
$L_{n+1}(\Z[\pi_1(M)])$, with $\T^{TOP}_{\partial}(M\times D^k)=[M \times D^k,\partial;G/TOP,*]$,
using the disjoint union $H$-space structure on $G/TOP$.
For a manifold with boundary $(M,\partial M)$ there are also rel and
rel $\partial$ versions. See Wall \cite[Theorem 10.8]{wall2}, Quinn
\cite{quinn1,quinn2}, Siebenmann \cite[Essay V, Appendix C]{ks},
Ranicki \cite{ranicki2,ranicki3}, Nicas \cite{nicas},
Cappell and Weinberger \cite{cw} and Kro \cite{kro}
for previous accounts of the topological structure set.

The {\it quadratic $\ST$-groups} $\ST_*(X)$ of a space $X$ were defined in
\cite[\S15]{ranicki4}, to fit into the algebraic surgery exact
sequence
$$\xymatrix@C-10pt{\dots \ar[r]& H_n(X;\LLdn) \ar[r]^-{\di{A}} &
 L_n(\Z[\pi_1(X)]) \ar[r]& \ST_n(X) \ar[r]&  H_{n-1}(X;\LLdn) \ar[r]& \dots}$$
with $A$ the assembly map. The functor
$$\{{\rm spaces}\} \to \{\Z\hbox{\rm -graded abelian groups}\}~;~
X \mapsto \ST_*(X)$$
is covariant and homotopy invariant: a map $f:X \to Y$ induces
morphisms $f_*:\ST_*(X) \to \ST_*(Y)$ depending only on the homotopy
class of $f$.  A homotopy equivalence $f$ induces isomorphisms $f_*$.

\begin{proposition} \label{global}
{\rm (Ranicki \cite[Theorem 18.5]{ranicki4})}
The algebraic and topological surgery exact sequences of an
$n$-dimensional manifold $M$ with $n
\geqslant 5$ are related by an isomorphism
$$\xymatrix@C-5pt{
\dots \ar[r]& L_{n+1}(\Z[\pi_1(M)])\ar@{=}[d] \ar[r]&
\ST^{TOP}(M) \ar[r]^-{\di{\eta}}\ar[d]^-{\di{s}}_{\di{\cong}}&
\T^{TOP}(M) \ar[r]^-{\di{\theta}}
\ar[d]^-{\di{t}}_{\di{\cong}} & L_n(\Z[\pi_1(M)]) \ar@{=}[d]\\
\dots \ar[r]& L_{n+1}(\Z[\pi_1(M)]) \ar[r]& \ST_{n+1}(M) \ar[r]&
 H_n(M;\LLdn) \ar[r]^-{\di{A}} & L_n(\Z[\pi_1(M)])}$$
 Similarly for a manifold with boundary.
 \end{proposition}
 \begin{proof} The isomorphism sending the topological normal
invariant to the algebraic normal invariant
$$\begin{array}{l}
t~=~[M]_{\LL}\cap -~:~\T^{TOP}(M)~=~[M,G/TOP]~=~H^0(M;\LLdn)
\xymatrix@C-5pt{\ar[r]^-{\di{\cong}}&} H_n(M;\LLdn)~;\\[1ex]
\hskip150pt \eta(f,b)\mapsto t(f,b)~=~[M]_{\LL}\cap \eta(f,b)
\end{array}$$
is the Poincar\'e duality isomorphism defined by cap product with the symmetric
$L$-theory fundamental class $[M]_{\LL} \in H_n(M;\LLup)$,
sending the topological normal invariant of a normal map $(f,b):N \to M$
to the algebraic normal invariant
$$t(f,b)~=~[M]_{\LL} \cap \eta(f,b) \in H_n(M;\LLdn)~.$$
Here, we are using the disjoint union $H$-space structure
on $G/TOP$ to define the abelian group structure, which corresponds
to the direct sum of quadratic Poincar\'e complexes in $\LLdn$.
The composite of $t$ and the quadratic $L$-theory assembly map $A$
$$At~:~\T^{TOP}(M)~\cong~H_n(M;\LLdn)\to L_n(\Z[\pi_1(M)])$$
sends a topological normal map $(f,b):N \to M$ to the surgery obstruction
$$At(f,b)~=~\sigma_*(f,b) \in L_n(\Z[\pi_1(M)])~.$$
The isomorphism
$$s~:~\ST^{TOP}(M) \xymatrix@C-5pt{\ar[r]^-{\di{\cong}}&} \ST_{n+1}(M)~;~
(N,f) \mapsto s(f)$$
sends a homotopy equivalence $f:N \to M$ to the
`manifold structure' $s(f)$ measuring the failure of the point
inverses $f^{-1}(x) \subset N$ ($x \in M$) to be acyclic.
(At this point it is worth recalling the $CE$ approximation theorem
of Siebenmann \cite{siebenmann}: for $n \geqslant 5$ a map of
$n$-dimensional topological manifolds $f:N \to M$ with contractible
point inverses is a homotopy equivalence which is homotopic to
a homeomorphism).
 \end{proof}

\begin{remark} {\rm Proposition \ref{global}
also applies in dimension $n=4$ in the case of a fundamental
group $\pi_1(M)$ which is good in the sense of Freedman and Quinn
\cite{fq}.}\hfill\qed
\end{remark}

\section{Composition formulae}

The composition formula for manifold structures is just
a local version of the composition formula of Ranicki \cite{ranicki3}
for surgery obstructions, so we start with that:

\begin{proposition}\label{c} {\rm (\cite[Proposition 4.3]{ranicki3})}
The composite of normal maps of $n$-dimensional geometric Poincar\'e complexes
$$(g,c)~:~P \to N~,~(f,b)~:~N \to M$$
is a normal map $(fg,bc):P \to M$ with surgery obstruction
$$\sigma_*(fg,bc)~=~\sigma_*(f,b)+f_*\sigma_*(g,c) \in L_n(\Z[\pi_1(M)])~.$$
\end{proposition}
\begin{proof}  Homotopy equivalent quadratic Poincar\'e complexes are cobordant.
The Umkehr chain maps are such that up to chain homotopy
$$(fg)^!~=~g^!f^!~:~C(\widetilde{M})
\xymatrix{\ar[r]^-{\displaystyle{f^!}} & }
C(\widetilde{N})
\xymatrix{\ar[r]^-{\displaystyle{g^!}} & } C(\widetilde{P})$$
and the quadratic kernels are such that up homotopy equivalence
$$(\C((fg)^!),\psi_{bc})~=~(\C(f^!),\psi_b)\oplus f_*(\C(g^!),\psi_c)~.$$
\end{proof}

\begin{remark} {\rm The surgery composition formula of Proposition
\ref{c} avoids the framing problems of the surgery composition formulae
for topological normal maps of Wall \cite[p.279]{wall2} and Jones
\cite[7.0]{jones}, being purely homotopy theoretic in nature.  See
Proposition \ref{match} below for the composition of topological normal
maps.}\hfill\qed
\end{remark}

The generalized $\LLdn$-homology groups $H_*(X;\LLdn)$ of a finite simplicial
complex $X$ were identified in Ranicki \cite[\S12]{ranicki4} with the quadratic
$L$-groups of an additive category $\A(\Z)_*(X)$ with chain duality in
which a quadratic Poincar\'e complex is a sheaf of quadratic Poincar\'e
complexes over $X$, i.e.  satisfies local Poincar\'e duality.
As in Ranicki and Weiss \cite{ranweiss} define a {\it $(\Z,X)$-module}
$A$ to be a $\Z$-module with a direct sum decomposition
$$A~=~\sum\limits_{\sigma \in X}A(\sigma)~,$$
with each $A(\sigma)$ based f.g.  free.  Let $\A(\Z)_*(X)$ be the additive category
with objects $(\Z,X)$-modules and morphisms $f:A \to B$ the abelian
group morphisms such that
$$f(A(\sigma)) \subseteq \sum\limits_{\tau\geqslant \sigma} B(\tau)~.$$
As in \cite[\S\S 5,13]{ranicki4} there is
defined a chain duality on $\A(\Z)_*(X)$ such that
$$H_n(X;\LLdn)~=~L_n(\A(\Z)_*(X))~.$$
An element $(C,\psi) \in H_n(X;\LLdn)$ is a
cobordism class of $n$-dimensional quadratic Poincar\'e complexes
$(C,\psi)$ in $\A(\Z)_*(X)$
with $C(\sigma)$ an $(n-\vert \sigma \vert)$-dimensional based f.g. free
$\Z$-module chain complex, subject to the additional requirement
(corresponding to the 1-connectivity of $\LLdn$)
that $C(\sigma)$ be contractible for $\sigma \in C$ with $\vert \sigma \vert=n$.
The quadratic structure is an element
$$\psi \in Q^X_n(C)~=~H_n(W\otimes_{\Z[\Z_2]}(C\otimes_{(\Z,X)}C))$$
with $C\otimes_{(\Z,X)}C \subseteq C \otimes_\Z C$ the $\Z[\Z_2]$-module
subcomplex defined by
$$C\otimes_{(\Z,X)}C~=~\sum\limits_{\sigma,\tau \in X,\sigma \cap \tau \neq \emptyset}
C(\sigma)\otimes_{\Z}C(\tau)~.$$
Let $\A(\Z[\pi_1(X)])$ be the additive category of based f.g.
free $\Z[\pi_1(X)]$-modules.
The assembly functor (\cite{ranweiss}, \cite[\S9]{ranicki4})
$$A~:~\A(\Z)_*(X) \to \A(\Z[\pi_1(X)])~;~B=\sum\limits_{\sigma \in X}B(\sigma)
\mapsto A(B)=\sum\limits_{\widetilde{\sigma}\in \widetilde{X}}B(p\widetilde{\sigma})$$
was defined using the universal covering projection $p:\widetilde{X}
\to X$.  An element $(C,\psi) \in H_n(X;\LLdn)$ is
a cobordism class of $n$-dimensional quadratic Poincar\'e complexes
$(C,\psi)$ in $\A(\Z)_*(X)$. An element $(C,\psi) \in \ST_{n+1}(X)$ is
a cobordism class of $n$-dimensional quadratic Poincar\'e complexes
$(C,\psi)$ in $\A(\Z)_*(X)$ such that the assembly based f.g.  free
$\Z[\pi_1(X)]$-module chain complex $A(C)$ is simple contractible.

We now recall from \cite{ranicki4} the construction of the algebraic
normal invariant $t(f,b) \in H_n(M;\LLdn)$ of a topological
normal map $(f,b):N \to M$ of $n$-dimensional manifolds.
By \cite[Corollary 17.6]{ranicki4} there exists a finite
simplicial complex $X$ with a homotopy
equivalence $h:M \to X$ such that both $h$ and $hf:N \to X$ are
topologically transverse across the dual cell decomposition
$\{D(\sigma,X)\vert \sigma \in X\}$ of the barycentric subdivision
$X'$.  The $\LLdn$-homology groups are homotopy invariant, so
$h_*:H_n(M;\LLdn)\to H_n(X;\LLdn)$ is an isomorphism.
The restrictions
$$(f(\sigma),b(\sigma))~=~(f,b)\vert~:~(hf)^{-1}D(\sigma,X) \to
h^{-1}D(\sigma,X)~~(\sigma \in X)$$
are normal maps of $(n-\vert \sigma \vert)$-dimensional geometric
Poincar\'e pairs (in fact manifolds) with boundaries the inverse images
of the boundaries of the dual cells
$$\partial D(\sigma,X)~=~\bigcup\limits_{\tau > \sigma}D(\tau,X)
\subset D(\sigma,X)~,$$
such that the degree 1 map is the assembly
$$f~=~\bigcup\limits_{\sigma \in X}f(\sigma)~:~
N~=~\bigcup\limits_{\sigma \in X}(hf)^{-1}D(\sigma,X)
\to M~=~\bigcup\limits_{\sigma \in X}h^{-1}D(\sigma,X)~.$$
The quadratic cycle $(C,\psi)$ has chain components
$$C(\sigma)~=~\C(f(\sigma)^!:C((hf)^{-1}D(\sigma,X))
\to C(h^{-1}D(\sigma,X)))~~(\sigma \in X)$$
with assembly the based f.g. free $\Z[\pi_1(X)]$-module chain
complex
$$A(C)~=~\C(\widetilde{f}^!:C(\widetilde{N}) \to C(\widetilde{M}))~,$$
and the quadratic components $\psi(\sigma)$ are the quadratic structures
of the normal maps $(f(\sigma),b(\sigma))$ used to define the
algebraic normal invariant
$$t(f,b) ~=~(C,\psi) \in H_n(M;\LLdn)~.$$
The algebraic manifold structure of a homotopy equivalence
$f:N \to M$ of $n$-dimensional manifolds is the cobordism class
$$s(f) ~=~(C,\psi) \in \ST_{n+1}(X)~=~\ST_{n+1}(M)$$
with $(C,\psi)$ (constructed as above) simple $\Z[\pi_1(M)]$-contractible by
virtue of $f$ being a homotopy equivalence. The bijection
of Proposition \ref{global}
$$\ST^{TOP}(M) \to \ST_{n+1}(M)~;~(N,f) \mapsto s(f)~=~(C,\psi)$$
sends the topological manifold structure to the algebraic manifold structure.

\begin{theorem}  \label{composite} Let $n\geqslant 5$.
The manifold structure of the composite $fg:P \to M$ of homotopy equivalences
$g:P \to N$, $f:N \to M$ of $n$-dimensional manifolds is given by
$$s(fg)~=~s(f)+f_*s(g) \in \ST^{TOP}(M)~=~\ST_{n+1}(M)~.$$
\end{theorem}
\begin{proof} Let $X$ be a finite simplicial complex with a homotopy
equivalence $h:M \to X$ such that  $h$ and $hf:N \to X$ and
$hfg:P \to X$ are
topologically transverse across the dual cell decomposition
$\{D(\sigma,X)\vert \sigma \in X\}$ of the barycentric subdivision
$X'$. Exactly as in the proof of Proposition \ref{c}
the chain complex kernels of the composite normal maps
$$\begin{array}{l}
((fg)(\sigma),(bc)(\sigma))~=~(f(\sigma),b(\sigma))(g(\sigma),c(\sigma))~:\\[1ex]
\hskip25pt
(hfg)^{-1}D(\sigma,X) \to (hf)^{-1}D(\sigma,X) \to h^{-1}D(\sigma,X)~~
(\sigma \in X)
\end{array}$$
split (up to $\Z$-module chain equivalence) as direct sums
$$\C((fg)(\sigma)^!)~=~\C(f(\sigma)^!) \oplus \C(g(\sigma)^!)~,$$
and similarly for the quadratic structures $\psi$. It follows that
$$s(fg)~=~s(f)+f_*s(g) \in \ST^{TOP}(M)~=~\ST_{n+1}(M)~=~\ST_{n+1}(X)$$
completing the proof.
\end{proof}

\begin{remark} {\rm
(i) Theorem \ref{composite} also applies in dimension $n=4$ in the case of a fundamental
group $\pi_1(M)$ which is good in the sense of Freedman and Quinn
\cite{fq}.\\
(ii) Theorem \ref{composite} is in fact true for all $n \geqslant
0$, provided that $f_*(P,g) \in \ST^{TOP}(M)$ is interpreted as an
element of $\ST_{n+1}(M)$ rather than as a manifold structure. The
algebraic manifold structure $s(f) \in \ST_{n+1}(M)$ is defined
for any homotopy equivalence $f:N \to M$ of $n$-dimensional
manifolds, and the proof of \ref{composite} gives
$s(fg)=s(f)+f_*s(g) \in \ST_{n+1}(M)$. \hfill\qed}
\end{remark}

For any degree 1 map $f:N \to M$ of $n$-dimensional geometric
Poincar\'e complexes there is defined a commutative square
$$\xymatrix{H^r(N) \ar@{<-}[r]^-{\di{f^*}}
\ar[d]_-{\di{[N]\cap -}}^-{\di{\cong}}
&H^r(M)\ar[d]^-{\di{f_*[N]\cap -}}_-{\di{\cong}} \\
H_{n-r}(N) \ar[r]^-{\di{f_*}} & H_{n-r}(M)}$$
with $f_*[N]=[M]\in H_n(M)$, so that the induced morphisms $f^*:H^r(N)
\to H^r(M)$ are split injections.  Similarly:

\begin{lemma} \label{diag}
For a topological normal map $(f,b):N \to M$ of $n$-dimensional manifolds
there is defined a commutative square
$$\xymatrix{[N,G/TOP]=H^0(N;\LLdn) \ar@{<-}[r]^-{\di{f^*}}
\ar[d]_-{\di{[N]_{\LL}\cap -}}^-{\di{\cong}}
&[M,G/TOP]=H^0(M;\LLdn)\ar[d]^-{\di{f_*[N]_{\LL}\cap -}}_-{\di{\cong}} \\
H_n(N;\LLdn) \ar[r]^-{\di{f_*}} & H_n(M;\LLdn)}$$
so that
$$f^*~:~[M,G/TOP] \to [N,G/TOP]~;~\eta \mapsto f^*\eta$$
is a split injection. \hfill\qed
\end{lemma}

\begin{remark} {\rm
(i) A degree 1 map of manifolds $f:N \to M$ does not in general
have degree 1 in $\LLup$-homology
$$f_*[N]_{\LL} \neq [M]_{\LL} \in H_n(M;\LLup)~.$$
(ii) If $(f,b):N \to M$ is a topological normal map then
$$f_*[N]_{\LL}~=~[M]_{\LL}+(1+T)t(f,b) \in H_n(M;\LLup)$$
is an $\LLdn$-theory fundamental class.}
\hfill\qed
\end{remark}

In general, it is not possible to compose topological normal maps.
However:

\begin{proposition} \label{match}
For topological normal maps $(f,b):N \to M$, $(g,c):P \to N$
such that
$$\eta(g,c) \in {\rm im}(f^*:[M,G/TOP] \to [N,G/TOP])$$
it is possible to define the composite
topological normal map $(fg,bc):P \to M$ with
$$bc~:~\nu_P \xymatrix{\ar[r]^-{\di{c}}&} \nu_N \oplus \nu_c
\xymatrix{\ar[r]^-{\di{b \oplus 1}}&} \nu_M \oplus \nu_b
\oplus f_!\nu_c~,$$
with normal invariant $f_!\eta(g,c) \in [M,G/TOP]$ the unique element such that
$$\eta(g,c)~=~f^*(f_!\eta(g,c)) \in [N,G/TOP]~.$$
{\rm (i)} The topological normal invariant of the composite is given by
$$\eta(fg,bc)~=~\eta(f,b) \oplus f_!\eta(g,c) \in [M,G/TOP]~.$$
{\rm (ii)} The algebraic normal invariant of the composite is given by
$$t(fg,bc)~=~t(f,b)+f_*t(g,c) \in H_n(M;\LLdn)~.$$
\end{proposition}
\begin{proof} (i) By construction.\\
(ii) Exactly as for the proof of Theorem \ref{composite}.
\end{proof}

\begin{remark} {\rm In particular, the hypothesis of Proposition \ref{match}
is satisfied in the special case when $f:N \to M$ is a homotopy equivalence,
with
$$f_!\eta(g,c)~=~(f^{-1})^*\eta(g,c) \in [M,G/TOP]~.$$
The composition formula in this case
$$\eta(fg,bc)~=~\eta(f,b) \oplus (f^{-1})^*\eta(g,c) \in [M,G/TOP]$$
was first obtained by Brumfiel \cite[Proposition 2.2]{brumfiel} and
Madsen, Taylor and Williams \cite[Lemma 2.5]{mtw}.}\hfill\qed
\end{remark}

In order to reconcile the composition formulae of Proposition \ref{match}
for the topological and algebraic normal invariants it is necessary
deal with the bijection
$$\eta~:~\T^{TOP}(M)~=~[M,G/TOP]\to H_n(M;\LLdn)~;~\eta(f,b)\mapsto t(f,b)~.$$
and compare the Whitney sum addition on the left with the
direct sum addition on the right, as will now be done in \S3.

\section{The Whitney sum and products}

It has long been known that the surgery obstruction function
$$\theta~:~\T^{TOP}(M)~=~[M,G/TOP] \to L_n(\Z[\pi_1(M)])$$
is a function between abelian groups which is not a morphism.
This was first observed for $PL$ normal invariants in the simply-connected
case by Cooke \cite[p.182]{cooke}, and in general by Wall
\cite[p.114]{wall2} :

\noindent
{\it We enter here a caveat to the reader.  In the situation of {\rm
(}10.6{\rm )}, we have the map
$$\theta~:~[X,G/PL] \rightarrow L_m\big(\pi(X)\big)$$
of abelian groups which satisfies, by definition, $\theta(0)=0$.
However, $\theta$ is {\sc NOT} in general a homomorphism.  The result
fails even in the closed, simply connected case with $4|m$, as one
readily sees by computing with Pontrjagin classes {\rm (}the simplest
example is the quaternion projective plane{\rm )}.}

\noindent
The classifying space $G/TOP$ has two distinct $H$-space structures,
one given by the Whitney sum of bundles and one given by the disjoint
union. The latter was originally defined by Sullivan using the
characteristic variety theorem and by Quinn using the
disjoint union of normal maps, and corresponds to the direct sum of quadratic
Poincar\'e complexes in $\LL_0 \simeq G/TOP$.  The abelian group
structure on $[M,G/TOP]$ defined by + matches the
abelian group structure on $H_n(M;\LLdn)$.
We shall now compare the composition formula for the structure
invariant (Theorem \ref{composite}) and the composition formula for the normal
invariant (Proposition \ref{match}), using the expression of the Whitney sum
$\oplus$ on $G/TOP$ as the addition
$$\oplus~:~\LL_0 \times \LL_0 \to \LL_0~;~(a,b) \mapsto a\oplus b=a+b+a\otimes b$$
(cf.  Ranicki \cite[p.295]{ranicki3}).

Recall the definition of the Whitney sum of bundles
$\eta:M \to BTOP(k)$, $\eta':M \to BTOP(k')$,
namely as the composite
$$\eta \oplus \eta'~:~M \xymatrix@C-5pt{\ar[r]^-{\di{\Delta}} &}
M \times M \xymatrix@C+5pt{\ar[r]^-{\di{\eta \times \eta'}} &}
\xymatrix@C-5pt{BTOP(k) \times BTOP(k') \ar[r] &} BTOP(k+k')~,$$
At this point, we need to recall the surgery product formula:

\begin{proposition} {\rm (Ranicki \cite[Proposition 8.1]{ranicki3})}\\
{\rm (i)} The product  $M \times M'$ of an $n$-dimensional geometric
Poincar\'e complex $M$ and an  $n'$-dimensional geometric Poincar\'e
complex $M'$ is an $(n+n')$-dimensional geometric Poincar\'e complex
$M \times M'$ with symmetric signature
$$\begin{array}{l}
\sigma^*(M \times M')~=~\sigma^*(M) \otimes_{\Z} \sigma^*(M')\\[1ex]
\hskip50pt
\in L^{n+n'}(\Z[\pi_1(M \times M')])~=~L^{n+n'}(\Z[\pi_1(M)]\otimes_{\Z}\Z[
\pi_1(M')])~.
\end{array}$$
{\rm (ii)} The product  of a normal map $(f,b):N \to M$ of
$n$-dimensional geometric Poincar\'e complexes and a normal map
$(f',b'):N' \to M'$ of $n'$-dimensional geometric Poincar\'e complexes
is a normal map $(f \times f',b \times b'):N \times N' \to M \times M'$
of $(n+n')$-dimensional geometric Poincar\'e complexes
with quadratic signature
$$\begin{array}{l}
\sigma_*(f\times f',b \times b')~=~\sigma_*(f,b)\otimes \sigma^*(M') +
\sigma^*(M) \otimes \sigma_*(f',b') +
\sigma_*(f,b)\otimes \sigma_*(f',b')\\[1ex]
\hskip50pt
\in L_{n+n'}(\Z[\pi_1(M \times M')])~=~L_{n+n'}(\Z[\pi_1(M)]\otimes_{\Z}\Z[
\pi_1(M')])~.
\end{array}$$
\end{proposition}
\begin{proof}  (i)
The symmetric complex of $M \times M'$ is such that up homotopy equivalence
$$(C(\widetilde{M \times M'}),\phi_{M \times M'})~=~
(C(\widetilde{M}),\phi_M)\otimes_{\Z}(C(\widetilde{M}'),\phi_{M'})$$
by the Eilenberg-Zilber theorem.
Homotopy equivalent symmetric Poincar\'e complexes are cobordant.\\
(ii) The quadratic kernel of $(f\times f',b \times b')$ is such that
up to homotopy equivalence
$$\begin{array}{l}
(\C((f\times f')^!),\psi_{b \times b'})~=\\[1ex]
((\C(f^!),\psi_b) \otimes_{\Z} (C(\widetilde{M}'),\phi_{M'})) \oplus
((C(\widetilde{M}),\phi_M)\otimes_{\Z}(\C({f'}^!),\psi_{b'}) )\\[1ex]
\hskip175pt \oplus ((\C(f^!),\psi_b) \otimes_{\Z}(\C({f'}^!),\psi_{b'}))~.
\end{array}$$
\end{proof}

Since $\LLup$ and $\LLdn$ are ring spectra and $\LLdn$ is an
$\LLup$-module spectrum there are defined cup and cap products for any
space $M$
$$\begin{array}{l}
\cup~:~H^p(M;\LLup) \otimes_\Z H^q(M;\LLup) \to H^{p+q}(M;\LLup)~,\\[1ex]
\cup~:~H^p(M;\LLdn) \otimes_\Z H^q(M;\LLdn) \to H^{p+q}(M;\LLdn)~,\\[1ex]
\cap~:~H_p(M;\LLup) \otimes_\Z H^q(M;\LLdn) \to H_{p-q}(M;\LLdn)~.
\end{array}$$
The symmetric $L$-theory fundamental class $[M]_{\LL} \in H_n(M;\LLup)$
of an $n$-dimensional manifold $M$ determines Poincar\'e duality
isomorphisms
$$[M]_{\LL} \cap - ~:~H^*(M;\LLdn) \to H_{n-*}(M;\LLdn)~;~
x \mapsto x_*~=~[M]_{\LL} \cap x$$
which can be used to define the intersection product
$$H_{n-p}(M;\LLdn) \otimes_\Z H_{n-q}(M;\LLdn) \to H_{n-p-q}(M;\LLdn)~;~
x_* \otimes y_* \mapsto x_*y_*=(x\cup y)_*$$
and similarly
$$H_{n-p}(M;\LLup) \otimes_\Z H_{n-q}(M;\LLdn) \to H_{n-p-q}(M;\LLdn)~.$$
We shall only be concerned with the intersection products in the special
case $p=q=0$
$$\begin{array}{l}
H_n(M;\LLdn) \otimes_\Z H_n(M;\LLdn) \to H_n(M;\LLdn)~,\\[1ex]
H_n(M;\LLup) \otimes_\Z H_n(M;\LLdn) \to H_n(M;\LLdn)~.
\end{array}$$
\begin{corollary} \label{product}
The product of an $n$-dimensional topological normal map
$(f,b):N \to M$ and an $n'$-dimensional topological normal map
$(f',b'):N' \to M'$ is an $(n+n')$-dimensional topological normal map
$$(f \times f',b \times b')~:~N \times N' \to M \times M'$$
with topological normal invariant
$$\eta(f \times f',b \times b')~=~
\eta(f,b)\times \eta(f',b')\in [M \times M',G/TOP]$$
and algebraic normal invariant
$$\begin{array}{ll}
t(f \times f',b \times b')&=~t(f,b)\otimes [M']_{\LL} +
[M]_{\LL} \otimes t(f',b') + t(f,b)\otimes t(f',b') \\[1ex]
&=~t(f,b)\otimes [M']_{\LL} + f_*[N]_\LL\otimes t(f',b')\in
H_{n+n'}(M \times M';\LLdn)
\end{array}$$
{\rm (}using $(1+T)t(f,b)=f_*[N]_{\LL}-[M]_{\LL}\in H_n(M;\LLup)${\rm )}.
\hfill\qed
\end{corollary}

The Whitney sum of normal invariants has the following
topological and algebraic properties:

\begin{theorem}\label{whitney}
Let $(f,b):N \to M$, $(f',b'):N' \to M$ be normal maps
of $n$-dimensional manifolds.\\
{\rm (i)} The Whitney sum of the topological normal invariants
$\eta(f,b),\eta(f',b') \in [M,G/TOP]$ is the topological
normal invariant
$$\eta(f,b) \oplus \eta(f',b')~=~\eta(f'',b'')
~=~\Delta^*(\eta(f,b) \times \eta(f',b'))  \in [M,G/TOP]$$
of the normal map $(f'',b'')$ obtained from the product
$$(f \times f',b \times b')~:~N \times N' \to M \times M$$
by transversality at the diagonal submanifold
$\Delta:M \subset M \times M$
$$(f'',b'')~=~(f \times f',b \times b')\vert~:~
N''~=~(f\times f')^{-1}(M) \to M~.$$
{\rm (ii)} The algebraic normal invariant of $(f'',b''):N'' \to M$ is
$$t(f'',b'')~=~t(f,b) + t(f',b') + t(f,b)t(f',b') \in H_n(M;\LLdn)$$
and surgery obstruction
$$\begin{array}{ll}
\sigma_*(f'',b'')&=~A(t(f'',b''))\\[1ex]
&=~\sigma_*(f,b)+\sigma_*(f',b')+A(t(f,b)t(f',b')) \in L_n(\Z[\pi_1(M)])~.
\end{array}$$
{\rm (iii)} Suppose that
$$\eta(f',b') \in {\rm im}(f^*:[N,G/TOP] \to [M,G/TOP])~,$$
so up to normal bordism there is a unique topological normal map
$(g,c):P \to N$ with
$$\begin{array}{l}
\eta(g,c)~=~f^*\eta(f',b') \in [N,G/TOP]~,\\[1ex]
t(g,c)~=~[N]_{\LL} \cap \eta(g,c) \in H_n(N;\LLdn)~,\\[1ex]
f_*t(g,c)~=~f_*[N]_{\LL} \cap t(f',b')~=~t(f',b')+t(f,b)t(f',b')
 \in H_n(M;\LLdn)~,\\[1ex]
\sigma_*(g,c)~=~A(f_*t(g,c))~=~\sigma_*(f',b')+A(t(f,b)t(f',b'))
\in L_n(\Z[\pi_1(M)])~.
\end{array}$$
Then
$$\eta(f'',b'')~=~\eta(f,b)\oplus \eta(f',b')~=~
\eta(fg,bc) \in [M,G/TOP]~,$$
so $(f'',b''):N'' \to M$ is normal bordant to $(fg,bc):P \to M$,
with
$$\begin{array}{l}
t(f'',b'')~=~t(f,b)+t(f',b')+t(f,b)t(f',b')\\[1ex]
\hphantom{t(f'',b'')~}=~
t(f,b)+f_*t(g,c)~=~t(fg,bc)~\in H_n(M;\LLdn)~,\\[1ex]
\sigma_*(f'',b'')~=~\sigma_*(f,b)+\sigma_*(f',b')+
A(t(f,b)t(f',b'))\\[1ex]
\hphantom{t(f'',b'')~}=~\sigma_*(f,b)+f_*\sigma_*(g,c)~=~
\sigma_*(fg,bc) \in L_n(\Z[\pi_1(M)])~.
\end{array}$$
\end{theorem}
\begin{proof} (i) By construction.\\
(ii) Pull back the product formulae of Corollary \ref{product}
along $\Delta:M \subset M \times M$, use $\nu_\Delta=\tau_M$,
$\nu_\Delta \oplus \nu_M \simeq *$ and the commutative square
$$\xymatrix@C+10pt{H^0(M\times M;\LLdn)
\ar[r]^-{\di[M \times M]_{\LL} \cap -}_-{\di{\cong}} \ar[d]^-{\di{\Delta^*}} &
H_{2n}(M\times M;\LLdn) \ar[d] \\
H^0(M;\LLdn) \ar[r]^-{\di[M]_{\LL} \cap -}_-{\di{\cong}}&
H_n(M;\LLdn) \cong H_{2n}(M\times M,M\times M \backslash \Delta(M);\LLdn)}$$\\
(iii) Immediate from (ii).\\
\end{proof}

We now reconcile the manifold structure composition formula of
Theorem \ref{composite} with the normal map composition formula
(Proposition \ref{match})
of Brumfiel \cite{brumfiel} and Madsen-Taylor-Williams \cite{mtw}:

\begin{corollary}\label{reconcile}
For homotopy equivalences of manifolds $f:N \to M$, $g:P \to N$
of $n$-dimensional manifolds
$$\eta(N,f)+\eta (f_*(P,g))~=~\eta(N,f)\oplus (f^{-1})^*\eta(P,g)
\in \T^{TOP}(M)~=~H_n(M;\LLdn)~.$$
\end{corollary}
\begin{proof}
Let $(f,b):N \to M$, $(g,c):P \to N$ be the corresponding
topological normal maps, with topological normal invariants
$$\eta(N,f)~=~\eta(f,b) \in \T^{TOP}(M)~,~
\eta(P,g)~=~\eta(g,c) \in \T^{TOP}(N)~.$$
If  $(f',b'):N' \to M$ is a topological normal map with
$$\eta(f',b')~=~(f^{-1})^*\eta(g,c) \in \T^{TOP}(M)~,$$
then by Proposition \ref{match}
the Whitney sum topological normal map of Theorem \ref{whitney} is
$$(f'',b'')~=~(fg,bc)~:~N''~=~P \to  M$$
up to normal bordism. Now
$$\sigma_*(f,b)~=~\sigma_*(g,c)~=~\sigma_*(fg,bc)~=~0 \in L_n(\Z[\pi_1(M)])~,$$
and Theorem \ref{whitney} (iii) gives
$$\begin{array}{l}
s(fg)~=~s(f)+f_*s(g) \in \ST_{n+1}(M)~,\\[1ex]
t(fg,bc)~=~t(f,b)+t(f',b')+t(f,b)t(f',b')\\[1ex]
\hphantom{t(fg,bc)~}=~t(f,b)+f_*t(g,c)\in H_n(M;\LLdn)~,\\[1ex]
\sigma_*(f',b')~=~-A(t(f,b)t(f',b')) \in L_n(\Z[\pi_1(M)])~.
\end{array}$$
Let $h:Q \to M$ be a homotopy equivalence of $n$-dimensional manifolds
with manifold structure $s(h)=f_*s(g) \in \ST_{n+1}(M)$.
The corresponding normal map $(h,d):Q \to M$ has algebraic normal invariant
$$\begin{array}{l}
t(h,d)~=~f_*t(g,c) ~=~ t(fg,bc)-t(f,b)~=~t(f',b')+t(f,b)t(f',b') \\[1ex]
\hskip25pt
\in {\rm im}(\ST_{n+1}(M) \to H_n(M;\LLdn))~=~
{\rm ker}(A:H_n(M;\LLdn) \to L_n(\Z[\pi_1(M)]))~.
\end{array}$$
The images of $(P,g) \in \ST^{TOP}(N)$ in the noncommutative diagram
$$\xymatrix@C+10pt{\ST^{TOP}(N)\ar[d]_-{\di{\eta}}
\ar[r]^-{\di{f_*}}_-{\di{\cong}} & \ST^{TOP}(M) \ar[d]^-{\di{\eta}} \\
\T^{TOP}(N) \ar[r]^-{\di{(f^{-1})^*}}_-{\di{\cong}} & \T^{TOP}(M)}$$
are given by
$$\begin{array}{l}
\eta (f_*(P,g))~=~\eta(Q,h)~=~t(h,d)~,\\[1ex]
(f^{-1})^*\eta(P,g)~=~(f^{-1})^*\eta(g,c)~=~t(f',b')
\in \T^{TOP}(M)~=~H_n(M;\LLdn)~,
\end{array}$$
differing by
$$t(h,d)-t(f',b')~=~t(f,b)t(f',b') \in \T^{TOP}(M)~=~H_n(M;\LLdn)~.$$
Thus
$$\begin{array}{ll}
\eta(N,f)+\eta (f_*(P,g))&=~t(f,b)+t(h,d)\\[1ex]
&=~t(f,b)+t(f',b')+t(f,b)t(f',b')\\[1ex]
&=~\eta(N,f)\oplus \eta(N',f')\\[1ex]
&=~\eta(N,f)\oplus (f^{-1})^*\eta(P,g)\in \T^{TOP}(M)~=~H_n(M;\LLdn)~.
\end{array}$$
\end{proof}

We conclude with a specific example, $M=S^p \times S^q$,
one of the two cases for which the manifold structure composition
formula $s(fg)=s(f)+f_*s(g)$ of Theorem \ref{composite}
is used by Kreck and L\"uck \cite{kl}.

\begin{example}
{\rm (i) Let $M=S^p\times S^q$ for $p,q \geqslant 2$,
so that $\pi_1(M)=\{1\}$. The assembly map in quadratic $L$-theory is given by
$$A~:~H_{p+q}(M;\LLdn)~=~L_p(\Z)\oplus L_q(\Z) \oplus L_{p+q}(\Z)
\to L_{p+q}(\Z)~;~(x,y,z) \mapsto z$$
and
$$\ST_{p+q+1}(M)~=~{\rm ker}(A)~=~L_p(\Z) \oplus L_q(\Z)~.$$
The addition and intersection pairing in $H_{p+q}(M;\LLdn)$ are given by
$$\begin{array}{l}
(x,y,z)+(x',y',z')~=~(x+x',y+y',z+z')~,\\[1ex]
(x,y,z)(x',y',z')~=~(0,0,xy'+x'y) \in H_{p+q}(M;\LLdn)~.
\end{array}$$
The product $L_p(\Z) \otimes L_q(\Z) \to L_{p+q}(\Z)$
factors through $L_p(\Z)\otimes L^q(\Z)$ and
the quadratic and symmetric $L$-groups of $\Z$ are given by
$$L_n(\Z)~=~\begin{cases} \Z \\ 0 \\ \Z_2 \\ 0
\end{cases}~,~L^n(\Z)~=~\begin{cases} \Z \\ \Z_2 \\ 0 \\ 0
\end{cases}~{\rm for}~n \equiv \begin{cases} 0 \\ 1 \\ 2\\ 3
\end{cases}~(\bmod\, 4)$$
so the intersection pairing
is non-zero only in the case $p\equiv q \equiv 0(\bmod\,4)$.
Given a topological normal map $(f,b):N \to M$ make $f$ transverse
regular at $S^p \times \{*\}$, $\{*\} \times S^q \subset M$
to obtain topological normal maps
$$\begin{array}{l}
(f_p,b_p)~=~(f,b)\vert~:~N_p~=~f^{-1}(S^p \times \{*\}) \to S^p ~,\\[1ex]
(f_q,b_q)~=~(f,b)\vert~:~N_q~=~f^{-1}(\{*\}\times S^q) \to S^q~.
\end{array}$$
and write the surgery obstructions as
$$(\sigma_*(f_p,b_p),\sigma_*(f_q,b_q),\sigma_*(f,b))~=~
(x_f,y_f,z_f) \in L_p(\Z) \oplus L_q(\Z) \oplus L_{p+q}(\Z)~.$$
The algebraic normal invariant defines a bijection
$$\T^{TOP}(M) \to L_p(\Z)\oplus L_q(\Z) \oplus L_{p+q}(\Z)~;~
\eta(f,b)\mapsto t(f,b)=(x_f,y_f,z_f)~.$$
The Whitney sum addition in $\T^{TOP}(M)$ corresponds to the addition
$$(x,y,z) \oplus (x',y',z')~=~(x+x',y+y',xy'+x'y+z+z') \in
L_p(\Z)\oplus L_q(\Z) \oplus L_{p+q}(\Z)~.$$
Given a homotopy equivalence $f:N \to M$ let
$(f,b):N \to M$ be the corresponding topological normal map.
The function
$$\ST^{TOP}(M) \to L_p(\Z)\oplus L_q(\Z)~;~s(f) \mapsto (x_f,y_f)$$
is a bijection, and
$$t(f,b)~=~(x_f,y_f,0) \in \T^{TOP}(M)~=~L_p(\Z)\oplus L_q(\Z) \oplus
L_{p+q}(\Z)~.$$
Given also a homotopy equivalence $g:P \to N$ with
corresponding topological normal map $(g,c):P \to N$ let
$$f_*s(g)~=~(x_g,y_g) \in \ST^{TOP}(M)~=~L_p(\Z)\oplus L_q(\Z)~,$$
so that
$$f_*t(g,c)~=~(x_g,y_g,0) \in \T^{TOP}(M)~=~L_p(\Z)\oplus L_q(\Z) \oplus
L_{p+q}(\Z)~.$$
As in the proof of
Corollary \ref{reconcile} let $(f',b'):N' \to M$ be a topological
normal map with topological normal invariant
$$\eta(f',b')~=~(f^{-1})^*\eta(g,c) \in \T^{TOP}(M)~,$$
let $h:Q \to N$ be a homotopy equivalence with
$$s(h)~=~f_*s(g)~=~(x_g,y_g) \in \ST^{TOP}(M)~=~\ST_{p+q+1}(M)~=~
L_p(\Z) \oplus L_q(\Z)~,$$
and let $(h,d):Q \to N$ be the corresponding topological normal map.
Then
$$\begin{array}{l}
\eta(f,b)\oplus \eta(f',b')~=~\eta(fg,bc) \in \T^{TOP}(M)~,\\[1ex]
t(f',b')~=~(x_g,y_g,-x_fy_g-x_gy_f)~,~t(h,d)~=~(x_g,y_g,0)~,\\[1ex]
t(fg,bc)~=~t(f,b)+f_*t(g,c)\\[1ex]
\hphantom{t(fg,bc)~}=~t(f,b)\oplus t(f',b')~=~t(f,b)+t(f',b')+t(f,b)t(f',b')\\[1ex]
\hphantom{t(fg,bc)~}=~(x_f+x_g,y_f+y_g,0)\in H_{p+q}(M;\LLdn)~
=~L_p(\Z) \oplus L_q(\Z) \oplus L_{p+q}(\Z)~,\\[1ex]
s(fg)~=~s(f)+f_*s(g)~=~(x_f+x_g,y_f+y_g) \in \ST^{TOP}(M)~=~
L_p(\Z) \oplus L_q(\Z)~.
\end{array}$$
(ii) For the simplest example of the non-additivity of the surgery obstruction
function
$$\theta~:~\T^{TOP}(M) \to L_n(\Z[\pi_1(M)])$$
with respect to $\oplus$ set $p=q=4$ in (i), and let
$f:N \to M=S^4 \times S^4$ be a homotopy equivalence with
$$\begin{array}{l}
s(f)~=~(x,y)\in \ST^{TOP}(M)~=~L_4(\Z) \oplus L_4(\Z)~,\\[1ex]
t(f,b)~=~(x,y,0)\in \T^{TOP}(M)~=~L_4(\Z) \oplus L_4(\Z)\oplus L_8(\Z)~,\\[1ex]
\theta(\eta(f,b))~=~\sigma_*(f,b)~=~A(t(f,b))~=~0\in L_8(\Z)
\end{array}$$
for any $x,y \in \Z \backslash \{0\}$. By (i)
a homotopy inverse $g=f^{-1}:P=M \to N$ is then such that
$$\begin{array}{l}
f_*s(g)~=~-s(f)~=~(-x,-y)\in \ST^{TOP}(M)~=~L_4(\Z) \oplus L_4(\Z)~,\\[1ex]
f_*t(g,c)~=~-t(f,b)~=~(-x,-y,0)\in \T^{TOP}(M)~=~L_4(\Z) \oplus L_4(\Z)\oplus L_8(\Z)
\end{array}$$
and a topological normal map $(f',b'):N' \to M$ with topological normal invariant
$\eta(f',b')=(f^{-1})^* \eta(g,c)$ has
$$\begin{array}{l}
t(f',b')~=~(-x,-y,2xy) \in \T^{TOP}(M)~=~L_4(\Z) \oplus L_4(\Z)\oplus L_8(\Z)~,\\[1ex]
\theta(\eta(f',b'))~=~\sigma_*(f',b')~=~A(t(f',b'))~=~2xy \neq 0
\in L_8(\Z)~=~\Z~.
\end{array}$$
The Whitney sum
$$\eta(f,b)\oplus \eta(f',b')~=~\eta(fg,bc)~=~\eta(1:M \to M)~=~0 \in \T^{TOP}(M)$$
has  surgery obstruction
$$\sigma_*(fg,bc)~=~\sigma_*(f,b)+\sigma_*(f',b')+A(t(f,b)t(f',b'))~=~0+2xy-2xy~=~0 \in L_8(\Z)~,$$
so
$$\theta(\eta(f,b)\oplus \eta(f',b'))=0 \neq \theta(\eta(f,b)) +
\theta(\eta(f',b'))=2xy \in L_8(\Z)=\Z~.$$
\hfill\qed
}
\end{example}


\begin{thebibliography}{99}

\bibitem{brumfiel} G.Brumfiel, {\it Homotopy equivalences of almost
smooth manifolds}, Comm. Math. Helv. 46, 381--407 (1971)

\bibitem{cw} S.Cappell and S.Weinberger, {\it A geometric
interpretation of Siebenmann's periodicity phenomenon}, Geometry
and topology (Athens, Ga., 1985), Lecture Notes in Pure and Appl.
Math. 105, 47--52, Dekker (1987)

\bibitem{cooke} G.E.Cooke, {\it The Hauptvermutung according to Casson
and Sullivan} (1968), published in The Hauptvermutung Book, 165--187,
Kluwer (1996)

\bibitem{fq} M.Freedman and F.Quinn, {\it Topology of 4-manifolds},
Princeton Mathematical Series 39, Princeton University Press (1990)

\bibitem{jones} L.Jones, {\it Patch spaces: a geometric represntation
for Poincar\'e spaces}, Ann. of Maths. 97, 306--343 (1973),
{\it Corrections for patch spaces}, ibid. 102, 183--185 (1975)

\bibitem{ks} R.Kirby and L.Siebenmann,
{\it Foundational Essays on Topological Manifolds, Smoothings, and
Triangulations}, Annals of Mathematics Studies 88, Princeton
(1977)

\bibitem{kl} M.Kreck and W.L\"uck, {\it Topological rigidity
for non-aspherical manifolds},
to appear in the Hirzebruch 80th birthday volume of the Quarterly
Journal of Pure and Applied Mathematics,
e-print arXiv:math.GT/0509238

\bibitem{kro} T.Kro, {\it Geometrically defined group structures in
3-dimensional surgery}, Oslo Cand. Scient. Thesis (2000),
http://www.maths.ntnu.no/$\widetilde{~}$toreak/hopp.pdf

\bibitem{mtw} I.Madsen, L.Taylor and B.Williams, {\it Tangential homotopy
equivalences}, Comm. Math. Helv. 55, 445-484 (1980)

\bibitem{nicas} A.Nicas, {\it Induction theorems for groups of homotopy manifold structures.}
Memoirs AMS 267 (1982)


\bibitem{quinn1} F.Quinn, {\it A geometric formulation of surgery},
Topology of manifolds, Proceedings 1969 Georgia Topology
Conference, Markham Press, 500--511 (1970)

\bibitem{quinn2} \bysame, {\it $B_{(TOP_n)^{\sim}}$ and the surgery
obstruction.} Bull.  A.M.S.  77, 596--600 (1971)


\bibitem{ranicki1} A.A.Ranicki, {\it The total surgery obstruction},
Proc. 1978 Arhus Topology Conference, Springer
Lecture Notes in Mathematics 763, 275--316 (1979)

\bibitem{ranicki2} \bysame, {\it The algebraic $L$-theory of surgery I.
Foundations}, Proc.  L.M.S.  (3) 40, 87--192 (1980)

\bibitem{ranicki3} \bysame, {\it The algebraic $L$-theory of surgery
II.  Applications to topology}, Proc.  L.M.S.  (3) 40, 193--283 (1980)


\bibitem{ranicki4} \bysame, {\it Algebraic $L$-theory and topological
manifolds}, Tracts in Mathematics 102, Cambridge  (1992)

\bibitem{ranweiss} \bysame\ and M.Weiss,
{\it Chain complexes and assembly}, Math. Z. 204, 157--185 (1990)

\bibitem{siebenmann} L.Siebenmann,
{\it Approximating cellular maps by homeomorphisms},
Topology 11, 271--294 (1972)

\bibitem{wall1} C.T.C.Wall, {\it Poincar\'e complexes},
Ann. Math. 86, 213--245 (1967)


\bibitem{wall2} \bysame, {\it Surgery on compact manifolds},
Academic Press (1970), 2nd edition, A.M.S. (1999)

\end{thebibliography}
\end{document}